\documentclass[paper=a4, USenglish, numbers=noenddot]{scrartcl} 
\usepackage{cmap}		
\usepackage[utf8]{inputenc}	
\usepackage[T1]{fontenc}	
\usepackage[USenglish]{babel}	
\usepackage[]{microtype}

\usepackage{amsmath}
\usepackage{amssymb}
\usepackage{graphicx}
\usepackage{mathtools}	

\usepackage{aliascnt}		
\usepackage{amsthm}

\usepackage{tikz}			
\usetikzlibrary{matrix,arrows,patterns,intersections,calc,decorations.pathmorphing}

\usepackage[pdftitle={The Veech group of the golden ladder},pdfauthor={Anja Randecker},pdfborder={0 0 0}]{hyperref}   
\usepackage[figure, table]{hypcap}		
\usepackage{pdfpages}

\usepackage{booktabs}
\usepackage{enumerate}
\usepackage{ragged2e, array}
\usepackage{xfrac}

\usepackage[disable]{todonotes}  
\usepackage[nameinlink]{cleveref}
\crefname{rem}{Remark}{Remarks}
\crefname{definition}{Definition}{Definitions}
\crefname{exa}{Example}{Examples}
\crefname{lem}{Lemma}{Lemmas}
\crefname{thm}{Theorem}{Theorems}
\crefname{thm1}{Theorem}{Theorems}
\crefname{prop}{Proposition}{Propositions}
\crefname{cor}{Corollary}{Corollaries}
\crefname{section}{Section}{Sections}
\crefname{subsection}{Subsection}{Subsections}
\crefname{figure}{Figure}{Figures}

\newtheoremstyle{break}
 {} 
 {} 
 {\itshape} 
 {} 
 {\bfseries} 
 {} 
 {\newline} 
 {\thmname{#1}\thmnumber{ #2}\thmnote{ (#3)}} 
 
\newtheoremstyle{breakdef}
 {} 
 {} 
 {} 
 {} 
 {\bfseries} 
 {} 
 {\newline} 
 {\thmname{#1}\thmnumber{ #2}\thmnote{ (#3)}} 
 
\newtheoremstyle{remark}
 {} 
 {} 
 {} 
 {} 
 {\itshape} 
 {.} 
 {0.5em} 
 {\thmname{#1}\thmnumber{ #2}\thmnote{ {\normalfont (}#3{\normalfont )}}} 

\theoremstyle{breakdef}

\theoremstyle{remark}

\newaliascnt{rem}{definition}  

\aliascntresetthe{rem}

\newaliascnt{exa}{definition}  

\aliascntresetthe{exa} 

\newaliascnt{cor}{definition}  

\aliascntresetthe{cor}  

\newaliascnt{lem}{definition}  
\newtheorem{lem}[lem]{Lemma}
\aliascntresetthe{lem}

\theoremstyle{break}

\newaliascnt{prop}{definition}  
\newtheorem{prop}[prop]{Proposition}
\aliascntresetthe{prop}

\newaliascnt{thm}{definition}  
\newtheorem{thm}[thm]{Theorem}
\aliascntresetthe{thm}

\DeclareMathOperator{\GL}{GL}
\DeclareMathOperator{\SL}{SL}
\DeclareMathOperator{\PSL}{PSL}

\renewcommand{\Re}{\mathrm{Re}}

\author{Anja Randecker}
\title{The Veech group of the golden ladder}
\date{\today}

\begin{document}

\maketitle

\begin{abstract}
 These are the notes of a talk that I gave at the Weihnachtsworkshop 2017 in Saarbrücken. It answered a question by Hooper and Treviño on the Veech group of the golden ladder, a translation surface of infinite type.
\end{abstract}

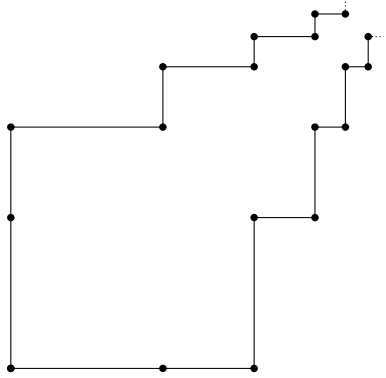
\begin{figure}
 \centering
 \begin{tikzpicture}[scale=2]
  \draw (0,0) -- (1,0)
   -- (1.6,0)
   -- (1.6,1)
   -- (2,1)
   -- (2,1.6)
   -- (2.2, 1.6)
   -- (2.2,2)
   -- (2.35,2)
   -- (2.35,2.2);
   \draw[densely dotted] (2.35,2.2) -- (2.45,2.2);
   
  \fill (0,0) circle (0.7pt);
  \fill (1,0) circle (0.7pt);
  \fill (1.6,0) circle (0.7pt);
  \fill (1.6,1) circle (0.7pt);
  \fill (2,1) circle (0.7pt);
  \fill (2,1.6) circle (0.7pt);
  \fill (2.2,1.6) circle (0.7pt);
  \fill (2.2,2) circle (0.7pt);
  \fill (2.35,2) circle (0.7pt);
  \fill (2.35,2.2) circle (0.7pt);

  \begin{scope}[rotate=90, yscale=-1]
   \draw (0,0) -- (1,0)
    -- (1.6,0)
    -- (1.6,1)
    -- (2,1)
    -- (2,1.6)
    -- (2.2, 1.6)
    -- (2.2,2)
    -- (2.35,2)
    -- (2.35,2.2);
   \draw[densely dotted] (2.35,2.2) -- (2.45,2.2);
    
   \fill (0,0) circle (0.7pt);
   \fill (1,0) circle (0.7pt);
   \fill (1.6,0) circle (0.7pt);
   \fill (1.6,1) circle (0.7pt);
   \fill (2,1) circle (0.7pt);
   \fill (2,1.6) circle (0.7pt);
   \fill (2.2,1.6) circle (0.7pt);
   \fill (2.2,2) circle (0.7pt);
   \fill (2.35,2) circle (0.7pt);
   \fill (2.35,2.2) circle (0.7pt);
  \end{scope}
 \end{tikzpicture}
 \caption{For the ladder surface, we identify opposite edges of this polygon.}
 \label{fig:ladder_surface_without_labels}
\end{figure}

A \emph{ladder surface} is a translation surface that we obtain when gluing opposite sides of an infinite-sided polygon as in \cref{fig:ladder_surface_without_labels}.

In the article \cite{hooper_trevino_17}, Hooper and Treviño consider a ladder surface with specific edge lengths that are related to the golden mean. They ask whether the symmetry group of this translation surface, the Veech group, is generated by a multi-twist and a reflection.

In these notes, we answer the question affirmatively, using a similar technique as for the baker's map surface (see \cite{chamanara_04, herrlich_randecker_16}).
In particular, the Veech group of the golden ladder is a non-elementary Fuchsian group of the second kind.

We show this result on the Veech group for a family of ladder surfaces where the golden ladder is a special case.

\section{Translation surfaces and Veech groups}

We first review the definitions that we use for these notes. In particular, note that in our definitions, translation surfaces are not necessarily of finite topological type, so we do not include the singularities into the surface.

A \emph{translation surface} $(X, \mathcal{A})$ is a connected two-dimensional manifold $X$ together with a translation structure $\mathcal{A}$ on $X$, that is, a maximal atlas on $X$ such that the transition functions are locally translations.
We can equip any translation surface $(X, \mathcal{A})$ with a flat metric by pulling back the Euclidean metric via the charts. The metric completion $\overline{X}$ may contain points which are not points of $X$ -- these are called \emph{singularities}.

There are three kinds of singularities:
A singularity is called \emph{cone angle singularity} if it has a punctured neighborhood in $X$ which is isometric to a finite translation covering of a once-punctured Euclidean disk. It is called \emph{infinite angle singularity} if it has a punctured neighborhood in $X$ which is isometric to an infinite translation covering of a once-punctured Euclidean disk. Any other singularity is called \emph{wild}.

\bigskip

A continuous map $f \colon (X,\mathcal{A}) \to (Y, \mathcal{B})$ of translation surfaces is called \emph{affine} if it is locally affine, that is, for every~$x\in X$ there exist charts $(U,\varphi)\in \mathcal{A}$ and $(V,\psi)\in \mathcal{B}$ with~$x\in U$ and $f(U)\subseteq V$ such that for every $z\in \varphi(U) \subseteq \mathbb{R}^2$, we have
\begin{equation*}
 \left( \psi\circ f \circ \varphi^{-1} \right)(z)= A\cdot z+t \text{ for a fixed } A\in \GL(2,\mathbb{R}) \text{ and a fixed } t\in\mathbb{R}^2 .
\end{equation*}

For an affine map $f$, the matrix $A$ is globally the same for all choices of charts. This matrix is called the \emph{derivative} of~$f$.
For a translation surface $(X,\mathcal{A})$, the \emph{Veech group}~$\GL^+(X,\mathcal{A})$ of $(X,\mathcal{A})$ is defined as the group of all derivatives of orientation-preserving homeomorphisms of $X$ that are affine with respect to $\mathcal{A}$.

For a translation surface with finite area, the elements of the affine group are area-preserving, hence the Veech group is a subgroup of $\SL(2,\mathbb{R})$.
To make use of fundamental domains of the action on $\mathbb{H}^2$ in \cref{sec:candidate}, we will subscribe to the point of view that the Veech group is in fact the image of $\GL^+(X, \mathcal{A})$ in $\PSL(2,\mathbb{R})$. This is also justified by the observation at the end of \cref{sec:Veech_group} that the Veech group of a ladder surface does not contain the element~$\begin{pmatrix}
 -1 & 0 \\
  0 & -1
\end{pmatrix}$.

\section{Ladder surfaces}

We describe a family of translation surfaces, called ladder surfaces. Note that the name comes from the fact that the gluing pattern resembles that of a one-sided infinite staircase, not from the topological type of a Jacob's ladder.

\begin{figure}
 \centering
 \begin{tikzpicture}[scale=2]
  \draw (0,0) -- node[below=-0.1]{$1$} (1,0)
   -- node[below=-0.1]{$\lambda$} (1.6,0)
   -- node[right=-0.1]{$1$} (1.6,1)
   -- node[below=-0.1]{$\lambda^2$} (2,1)
   -- node[right=-0.1]{$\lambda$} (2,1.6)
   -- (2.2, 1.6) 
   -- node[right=-0.1]{$\lambda^2$} (2.2,2)
   -- (2.35,2)
   -- (2.35,2.2);
   \draw[densely dotted] (2.35,2.2) -- (2.45,2.2);
   
  \fill (0,0) circle (0.7pt);
  \fill (1,0) circle (0.7pt);
  \fill (1.6,0) circle (0.7pt);
  \fill (1.6,1) circle (0.7pt);
  \fill (2,1) circle (0.7pt);
  \fill (2,1.6) circle (0.7pt);
  \fill (2.2,1.6) circle (0.7pt);
  \fill (2.2,2) circle (0.7pt);
  \fill (2.35,2) circle (0.7pt);
  \fill (2.35,2.2) circle (0.7pt);
   
  \draw (2.1, 1.58) -- (2.4, 1.4) node[right=-0.1]{$\lambda^3$};

  \begin{scope}[rotate=90, yscale=-1]
   \draw (0,0) -- node[left=-0.1]{$1$} (1,0)
    -- node[left=-0.1]{$\lambda$} (1.6,0)
    -- node[above=-0.1]{$1$} (1.6,1)
    -- node[left=-0.1]{$\lambda^2$} (2,1)
    -- node[above=-0.1]{$\lambda$} (2,1.6)
    -- (2.2, 1.6) 
    -- node[above=-0.1]{$\lambda^2$} (2.2,2)
    -- (2.35,2)
    -- (2.35,2.2);
   \draw[densely dotted] (2.35,2.2) -- (2.45,2.2);
    
   \fill (0,0) circle (0.7pt);
   \fill (1,0) circle (0.7pt);
   \fill (1.6,0) circle (0.7pt);
   \fill (1.6,1) circle (0.7pt);
   \fill (2,1) circle (0.7pt);
   \fill (2,1.6) circle (0.7pt);
   \fill (2.2,1.6) circle (0.7pt);
   \fill (2.2,2) circle (0.7pt);
   \fill (2.35,2) circle (0.7pt);
   \fill (2.35,2.2) circle (0.7pt);
    
   \draw (2.1, 1.58) -- (2.3, 1.3) node[above=-0.1]{$\lambda^3$};
  \end{scope}
 \end{tikzpicture}
 \caption{For the ladder surface with parameter $\lambda$, we identify edges that are parallel and have the same length (labels are indicating length).}
 \label{fig:ladder_surface}
\end{figure}
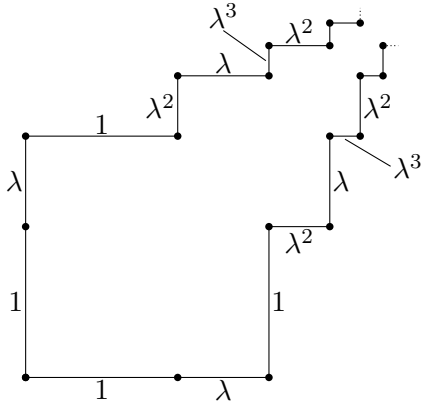

We define a translation surface for every $\lambda \in (0,\infty)$ with a gluing pattern and edge lengths as in \autoref{fig:ladder_surface}.
For $\lambda = 0$, the surface is a torus; for $\lambda = 1$, it is a one-sided version of the infinite staircase as in \cite{hooper_hubert_weiss_13, malaga-sabogal_troubetzkoy_19}.

We are particularly interested in the case of $\lambda \in (0,1)$. Then the area of the surface is
\begin{equation*}
 1^2 + 2 \cdot \lambda \cdot 1 + \lambda^2 + 2 \cdot \lambda^2 \cdot \lambda + \ldots
 = \sum_{k=0}^\infty \lambda^{2k} + 2 \cdot \sum_{k=0}^\infty \lambda^{2k+1}
 < \infty
\end{equation*}
and it has exactly one singularity, which is wild.
In the case of $\lambda \geq 1$, the ladder surface has infinite area and three infinite angle singularities.
In all cases, the surface has one end and infinite genus. Hence, the topological type is that of a Loch Ness~monster.

A ladder surface has many symmetries if its parameter $\lambda$ fulfills the following condition:
Let $k, l \in \mathbb{Z}$ be coprime numbers with $k > l > 0$. Then we define $\lambda_{k,l}$ to be the greater (i.e.\ the positive) solution of the equation $k \cdot (\lambda +1) = l \cdot (\lambda^{-1} + 1 + \lambda)$.
In particular for~$l = 1$ and $k = 2$, $\lambda_{2,1}$ is the inverse of the golden mean. This is the case of the golden ladder that was considered in \cite{hooper_trevino_17}.

Ladder surfaces arise in different ways from previously known constructions and examples.
Hooper and Treviño describe the golden ladder as a translation surface associated to an infinite grid graph construction.

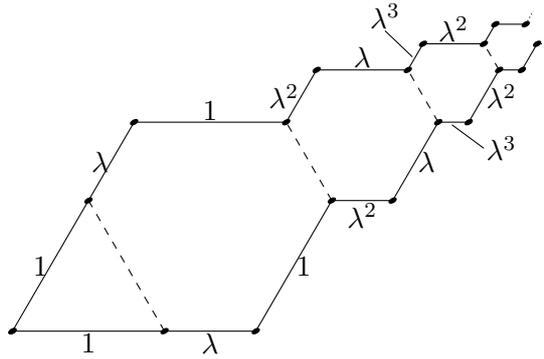
\begin{figure}
 \centering
 \begin{tikzpicture}[scale=2, cm={1,0,0.5,.866,(0,0)}]
  \draw (0,0) -- node[below=-0.1]{$1$} (1,0)
   -- node[below=-0.1]{$\lambda$} (1.6,0)
   -- node[right=-0.1]{$1$} (1.6,1)
   -- node[below=-0.1]{$\lambda^2$} (2,1)
   -- node[right=-0.1]{$\lambda$} (2,1.6)
   -- (2.2, 1.6) 
   -- node[right=-0.1]{$\lambda^2$} (2.2,2)
   -- (2.35,2)
   -- (2.35,2.2);
   \draw[densely dotted] (2.35,2.2) -- (2.45,2.2);
   
  \fill (0,0) circle (0.7pt);
  \fill (1,0) circle (0.7pt);
  \fill (1.6,0) circle (0.7pt);
  \fill (1.6,1) circle (0.7pt);
  \fill (2,1) circle (0.7pt);
  \fill (2,1.6) circle (0.7pt);
  \fill (2.2,1.6) circle (0.7pt);
  \fill (2.2,2) circle (0.7pt);
  \fill (2.35,2) circle (0.7pt);
  \fill (2.35,2.2) circle (0.7pt);
   
  \draw (2.1, 1.58) -- (2.4, 1.4) node[right=-0.1]{$\lambda^3$};

  \begin{scope}[rotate=90, yscale=-1]
   \draw (0,0) -- node[left=-0.1]{$1$} (1,0)
    -- node[left=-0.1]{$\lambda$} (1.6,0)
    -- node[above=-0.1]{$1$} (1.6,1)
    -- node[left=-0.1]{$\lambda^2$} (2,1)
    -- node[above=-0.1]{$\lambda$} (2,1.6)
    -- (2.2, 1.6) 
    -- node[above=-0.1]{$\lambda^2$} (2.2,2)
    -- (2.35,2)
    -- (2.35,2.2);
   \draw[densely dotted] (2.35,2.2) -- (2.45,2.2);
    
   \fill (0,0) circle (0.7pt);
   \fill (1,0) circle (0.7pt);
   \fill (1.6,0) circle (0.7pt);
   \fill (1.6,1) circle (0.7pt);
   \fill (2,1) circle (0.7pt);
   \fill (2,1.6) circle (0.7pt);
   \fill (2.2,1.6) circle (0.7pt);
   \fill (2.2,2) circle (0.7pt);
   \fill (2.35,2) circle (0.7pt);
   \fill (2.35,2.2) circle (0.7pt);
    
   \draw (2.1, 1.58) -- (2.3, 1.3) node[above=-0.1]{$\lambda^3$};
  \end{scope}
  
  \draw[dashed] (1,0) -- (0,1);
  \draw[dashed] (1.6,1) -- (1,1.6);
  \draw[dashed] (2,1.6) -- (1.6,2);
  \draw[dashed] (2.2,2) -- (2,2.2);
 \end{tikzpicture}
 \caption{The ladder surface with parameter $\lambda$ after applying an affine transformation.}
 \label{fig:bouw-moeller_surface}
\end{figure}

We can also apply the matrix
$\begin{pmatrix}
  1 & \sfrac{1}{2} \\
  0 & \sfrac{\sqrt{3}}{2}
 \end{pmatrix}$
(a shear combined with a scaling) and obtain a translation surface as in \autoref{fig:bouw-moeller_surface}.
This translation surface can be thought of as being glued from semi-regular hexagons. Here, semiregular means that every second edge has the same length. The first hexagon actually has the edge lengths $1$ and $0$, alternatingly.

Finite translation surfaces that are glued together from semiregular polygons have been described by Hooper in \cite{hooper_13c}, based on earlier work of Bouw and Möller in \cite{bouw_moeller_10} where the translation surfaces are given by a different description.

It is convenient to switch between the description as a ladder surface and as gluing of semi-regular hexagons to determine the Veech group.
In particular, we can see from the description by semi-regular hexagon that there exists an elliptic element of order $3$ in the Veech group.

\section{Three parabolic elements of the Veech group} \label{sec_cylinder_decompositions}

For a translation surface $(X,\mathcal{A})$ with at least one singularity, a \emph{cylinder} in $(X,\mathcal{A})$ of \emph{circumference} $w>0$ and \emph{height} $h>0$ is an open subset of $X$ which is isometric to a Euclidean cylinder $\mathbb{R}/w\mathbb{Z} \times (0,h)$.
The \emph{modulus} of a cylinder is the ratio $\frac{w}{h}$ of circumference and height.

A \emph{saddle connection} of $(X,\mathcal{A})$ is a geodesic segment in $\overline{X}$ from one singularity to another or the same singularity which does not contain singularities in its interior. 
If a cylinder can be extended to a maximal cylinder then this maximal cylinder is bounded by saddle connections.
A \emph{cylinder decomposition} of $(X,\mathcal{A})$ is a collection of maximal cylinders in $(X,\mathcal{A})$ such that the closures of the cylinders in $X$ cover $X$ and such that each two cylinders are disjoint.
The \emph{direction} of the cylinder decomposition is the direction of the corresponding bounding saddle connections.

We will use the following standard result to find elements of the Veech group. It was first proven in \cite{veech_89} as an equivalence for translation surfaces of finite type. However, the proof of the implication that is stated here literally works for the general~case.

\begin{prop}[Cylinder decompositions and parabolic elements {[Veech]}] \label{prop_cylinder_decomposition_parabolic_element}
 Let $(X,\mathcal{A})$ be a translation surface and $\left\{z_n\right\}$ a cylinder decomposition such that the cylinder~$z_n$ has height $h_n$ and circumference $w_n$. If the inverse moduli $\frac{h_n}{w_n}$ are commensurable, that is, if there exists an $m\in \mathbb{R}$ such that each inverse modulus is an integer multiple of~$m$, then the Veech group contains a parabolic element conjugated to the matrix
 $\begin{pmatrix}
  1 & \sfrac{1}{m} \\
  0 & 1  
 \end{pmatrix}$.
\end{prop}

\begin{figure}
 \centering
 \begin{tikzpicture}[scale=2]
  \draw[pattern color=gray!60!white, pattern=north east lines] (0,0) -- (1.6,0) -- (1.6,1) -- (0,1) -- (0,0);
  \draw[pattern color=gray!60!white, pattern=dots] (0,1) -- (2,1) -- (2,1.6) -- (0,1.6) -- (0,1);
  \draw[pattern color=gray!60!white, pattern=bricks] (1,1.6) -- (2.2,1.6) -- (2.2,2) -- (1,2) -- (1,1.6);
  \draw[pattern color=gray!60!white, pattern=horizontal lines] (1.6,2) -- (2.35,2) -- (2.35,2.2) -- (1.6,2.2) -- (1.6,2); 
 
  \draw (0,0) -- node[below=-0.1]{$1$} (1,0)
   -- node[below=-0.1]{$\lambda$} (1.6,0)
   -- node[right=-0.1]{$1$} (1.6,1)
   -- node[below=-0.1]{$\lambda^2$} (2,1)
   -- node[right=-0.1]{$\lambda$} (2,1.6)
   -- (2.2, 1.6) 
   -- node[right=-0.1]{$\lambda^2$} (2.2,2)
   -- (2.35,2)
   -- (2.35,2.2);
   \draw[densely dotted] (2.35,2.2) -- (2.45,2.2);
   
  \fill (0,0) circle (0.7pt);
  \fill (1,0) circle (0.7pt);
  \fill (1.6,0) circle (0.7pt);
  \fill (1.6,1) circle (0.7pt);
  \fill (2,1) circle (0.7pt);
  \fill (2,1.6) circle (0.7pt);
  \fill (2.2,1.6) circle (0.7pt);
  \fill (2.2,2) circle (0.7pt);
  \fill (2.35,2) circle (0.7pt);
  \fill (2.35,2.2) circle (0.7pt);
   
  \draw (2.1, 1.58) -- (2.4, 1.4) node[right=-0.1]{$\lambda^3$};

  \begin{scope}[rotate=90, yscale=-1]
   \draw (0,0) -- node[left=-0.1]{$1$} (1,0)
    -- node[left=-0.1]{$\lambda$} (1.6,0)
    -- node[above=-0.1]{$1$} (1.6,1)
    -- node[left=-0.1]{$\lambda^2$} (2,1)
    -- node[above=-0.1]{$\lambda$} (2,1.6)
    -- (2.2, 1.6) 
    -- node[above=-0.1]{$\lambda^2$} (2.2,2)
    -- (2.35,2)
    -- (2.35,2.2);
   \draw[densely dotted] (2.35,2.2) -- (2.45,2.2);
    
   \fill (0,0) circle (0.7pt);
   \fill (1,0) circle (0.7pt);
   \fill (1.6,0) circle (0.7pt);
   \fill (1.6,1) circle (0.7pt);
   \fill (2,1) circle (0.7pt);
   \fill (2,1.6) circle (0.7pt);
   \fill (2.2,1.6) circle (0.7pt);
   \fill (2.2,2) circle (0.7pt);
   \fill (2.35,2) circle (0.7pt);
   \fill (2.35,2.2) circle (0.7pt);
    
   \draw (2.1, 1.58) -- (2.3, 1.3) node[above=-0.1]{$\lambda^3$};
  \end{scope}
 \end{tikzpicture}
 \label{fig:horizontal_cylinder_decomposition}
 \caption{Ladder surface with a horizontal cylinder decomposition.}
\end{figure}
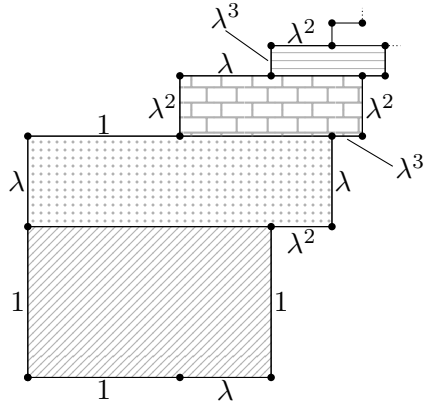

\autoref{fig:horizontal_cylinder_decomposition} shows a horizontal cylinder decomposition on the ladder surface. It consists of infinitely many cylinders where all except the bottom one are similar.
The bottom cylinder has a height of $1$ and a circumference of $1+\lambda$, hence the modulus is $1+\lambda$.
The cylinder above has a height of $\lambda$ and a circumference of $1+\lambda+\lambda^2$, hence the modulus is~$\lambda^{-1}+1+\lambda$. For $\lambda = \lambda_{k,l}$, we get that the moduli of all the cylinders but the bottom one are $\frac{k}{l} \cdot (1 + \lambda)$.
Then \autoref{prop_cylinder_decomposition_parabolic_element} implies that there exists a parabolic element in the Veech group which preserves the horizontal direction and acts as a (multiple) Dehn twist on every cylinder in the horizontal cylinder decomposition.
The derivative of it is
$\begin{pmatrix}
  1 & k(1+\lambda) \\
  0 & 1  
\end{pmatrix}$ and it twists the bottom cylinder $k$ times and all other cylinders $l$ times.

\begin{figure}[bhtp]
 \centering
 \begin{tikzpicture}[scale=1.4, cm={1,0,0.5,.866,(0,0)}]
  \draw[pattern color=gray!60!white, pattern=north east lines] (0,0) -- (1.6,0) -- (1.6,1) -- (0,1) -- (0,0);
  \draw[pattern color=gray!60!white, pattern=dots] (0,1) -- (2,1) -- (2,1.6) -- (0,1.6) -- (0,1);
  \draw[pattern color=gray!60!white, pattern=bricks] (1,1.6) -- (2.2,1.6) -- (2.2,2) -- (1,2) -- (1,1.6);
  \draw[pattern color=gray!60!white, pattern=horizontal lines] (1.6,2) -- (2.35,2) -- (2.35,2.2) -- (1.6,2.2) -- (1.6,2); 
 
  \draw (0,0) -- (1,0)
   -- (1.6,0)
   -- (1.6,1)
   -- (2,1)
   -- (2,1.6)
   -- (2.2, 1.6)
   -- (2.2,2)
   -- (2.35,2)
   -- (2.35,2.2);
   \draw[densely dotted] (2.35,2.2) -- (2.45,2.2);
   
  \fill (0,0) circle (0.7pt);
  \fill (1,0) circle (0.7pt);
  \fill (1.6,0) circle (0.7pt);
  \fill (1.6,1) circle (0.7pt);
  \fill (2,1) circle (0.7pt);
  \fill (2,1.6) circle (0.7pt);
  \fill (2.2,1.6) circle (0.7pt);
  \fill (2.2,2) circle (0.7pt);
  \fill (2.35,2) circle (0.7pt);
  \fill (2.35,2.2) circle (0.7pt);
   
  \begin{scope}[rotate=90, yscale=-1]
   \draw (0,0) -- (1,0)
    -- (1.6,0)
    -- (1.6,1)
    -- (2,1)
    -- (2,1.6)
    -- (2.2, 1.6)
    -- (2.2,2)
    -- (2.35,2)
    -- (2.35,2.2);
   \draw[densely dotted] (2.35,2.2) -- (2.45,2.2);
    
   \fill (0,0) circle (0.7pt);
   \fill (1,0) circle (0.7pt);
   \fill (1.6,0) circle (0.7pt);
   \fill (1.6,1) circle (0.7pt);
   \fill (2,1) circle (0.7pt);
   \fill (2,1.6) circle (0.7pt);
   \fill (2.2,1.6) circle (0.7pt);
   \fill (2.2,2) circle (0.7pt);
   \fill (2.35,2) circle (0.7pt);
   \fill (2.35,2.2) circle (0.7pt);
  \end{scope}
  
  \begin{scope}[xshift=3.3cm] 
  \draw[pattern color=gray!60!white, pattern=north east lines] (0,0) -- (1,0) -- (1,1.6) -- (0,1.6) -- (0,0);
  \draw[pattern color=gray!60!white, pattern=dots] (1,0) -- (1.6,0) -- (1.6,2) -- (1,2) -- (1,0);
  \draw[pattern color=gray!60!white, pattern=bricks] (1.6,1) -- (2,1) -- (2,2.2) -- (1.6,2.2) -- (1.6,1);
  \draw[pattern color=gray!60!white, pattern=horizontal lines] (2,1.6) -- (2.2,1.6) -- (2.2,2.35) -- (2,2.35) -- (2,1.6);
 
  \draw (0,0) -- (1,0)
   -- (1.6,0)
   -- (1.6,1)
   -- (2,1)
   -- (2,1.6)
   -- (2.2, 1.6)
   -- (2.2,2)
   -- (2.35,2)
   -- (2.35,2.2);
   \draw[densely dotted] (2.35,2.2) -- (2.45,2.2);
   
  \fill (0,0) circle (0.7pt);
  \fill (1,0) circle (0.7pt);
  \fill (1.6,0) circle (0.7pt);
  \fill (1.6,1) circle (0.7pt);
  \fill (2,1) circle (0.7pt);
  \fill (2,1.6) circle (0.7pt);
  \fill (2.2,1.6) circle (0.7pt);
  \fill (2.2,2) circle (0.7pt);
  \fill (2.35,2) circle (0.7pt);
  \fill (2.35,2.2) circle (0.7pt);
   
  \begin{scope}[rotate=90, yscale=-1]
   \draw (0,0) -- (1,0)
    -- (1.6,0)
    -- (1.6,1)
    -- (2,1)
    -- (2,1.6)
    -- (2.2, 1.6)
    -- (2.2,2)
    -- (2.35,2)
    -- (2.35,2.2);
   \draw[densely dotted] (2.35,2.2) -- (2.45,2.2);
    
   \fill (0,0) circle (0.7pt);
   \fill (1,0) circle (0.7pt);
   \fill (1.6,0) circle (0.7pt);
   \fill (1.6,1) circle (0.7pt);
   \fill (2,1) circle (0.7pt);
   \fill (2,1.6) circle (0.7pt);
   \fill (2.2,1.6) circle (0.7pt);
   \fill (2.2,2) circle (0.7pt);
   \fill (2.35,2) circle (0.7pt);
   \fill (2.35,2.2) circle (0.7pt);
  \end{scope}
  \end{scope}
  
  \begin{scope}[xshift=6.6cm] 
  \draw[pattern color=gray!60!white, pattern=north east lines] (0,0) -- (1,0) -- (0,1) -- (0,0);
  \draw[pattern color=gray!60!white, pattern=north east lines] (1.6,0) -- (1.6,1) -- (1,1.6) -- (0,1.6) -- (1.6,0);
  \draw[pattern color=gray!60!white, pattern=dots] (1,0) -- (1.6,0) -- (0,1.6) -- (0,1) -- (1,0);
  \draw[pattern color=gray!60!white, pattern=dots] (2,1) -- (2,1.6) -- (1.6,2) -- (1,2) -- (2,1);
  \draw[pattern color=gray!60!white, pattern=bricks] (1.6,1) -- (2,1) -- (1,2) -- (1,1.6) -- (1.6,1);
  \draw[pattern color=gray!60!white, pattern=bricks] (2.2,1.6) -- (2.2,2) -- (2,2.2) -- (1.6,2.2) -- (2.2,1.6);
  \draw[pattern color=gray!60!white, pattern=horizontal lines] (2,1.6) -- (2.2,1.6) -- (1.6,2.2) -- (1.6,2) -- (2,1.6);
  \draw[pattern color=gray!60!white, pattern=horizontal lines] (2.35,2) -- (2.35,2.2) -- (2.2,2.35) -- (2,2.35) -- (2.35,2); 
 
  \draw (0,0) -- (1,0)
   -- (1.6,0)
   -- (1.6,1)
   -- (2,1)
   -- (2,1.6)
   -- (2.2, 1.6)
   -- (2.2,2)
   -- (2.35,2)
   -- (2.35,2.2);
   \draw[densely dotted] (2.35,2.2) -- (2.45,2.2);
   
  \fill (0,0) circle (0.7pt);
  \fill (1,0) circle (0.7pt);
  \fill (1.6,0) circle (0.7pt);
  \fill (1.6,1) circle (0.7pt);
  \fill (2,1) circle (0.7pt);
  \fill (2,1.6) circle (0.7pt);
  \fill (2.2,1.6) circle (0.7pt);
  \fill (2.2,2) circle (0.7pt);
  \fill (2.35,2) circle (0.7pt);
  \fill (2.35,2.2) circle (0.7pt);
   
  \begin{scope}[rotate=90, yscale=-1]
   \draw (0,0) -- (1,0)
    -- (1.6,0)
    -- (1.6,1)
    -- (2,1)
    -- (2,1.6)
    -- (2.2, 1.6)
    -- (2.2,2)
    -- (2.35,2)
    -- (2.35,2.2);
   \draw[densely dotted] (2.35,2.2) -- (2.45,2.2);
    
   \fill (0,0) circle (0.7pt);
   \fill (1,0) circle (0.7pt);
   \fill (1.6,0) circle (0.7pt);
   \fill (1.6,1) circle (0.7pt);
   \fill (2,1) circle (0.7pt);
   \fill (2,1.6) circle (0.7pt);
   \fill (2.2,1.6) circle (0.7pt);
   \fill (2.2,2) circle (0.7pt);
   \fill (2.35,2) circle (0.7pt);
   \fill (2.35,2.2) circle (0.7pt);
  \end{scope}
  \end{scope}
 \end{tikzpicture}
 \caption{Three cylinder decompositions in the description as semiregular hexagons.}
 \label{fig:bouw-moeller_surface_cylinders}
\end{figure}
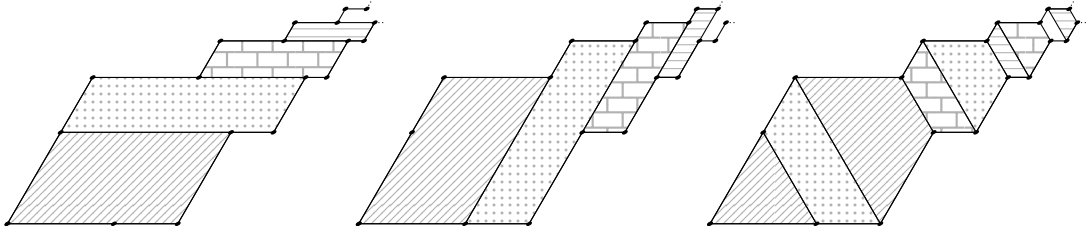

Analogously, we have another parabolic element in the Veech group, coming from a vertical cylinder decomposition.
When considering the description by semiregular hexagons, we can even see three cylinder decompositions directly as in \cref{fig:bouw-moeller_surface_cylinders}: The horizontal decomposition in this description corresponds to the horizontal decomposition we described before. Using the elliptic element of order $3$, we also have two diagonal directions of cylinder decompositions which correspond to the vertical decomposition and a decomposition of slope $-1$ in the description as ladder surface.

\bigskip

As a converse to the existence of these three cylinder directions, we prove that a large set of directions cannot occur as eigen directions of parabolic elements in the Veech group. Note that there can still exist cylinder decompositions in these directions.

\begin{lem}\label{lem_forbidden_eigen_directions}
 There is no parabolic element in the Veech group of the ladder surface with parameter $\lambda \in (0,1)$ whose eigen direction corresponds to a slope in the interval $(\lambda, \lambda^{-1})$.
 
\begin{proof}
 Consider the point $S = \left(\frac{1}{1-\lambda}, \frac{1}{1-\lambda} \right)$ in the metric completion of the ladder surface. As $\frac{1}{1-\lambda} = \sum_{n=0}^\infty \lambda^n$, this point sits in the top right corner of the ladder surface.
  
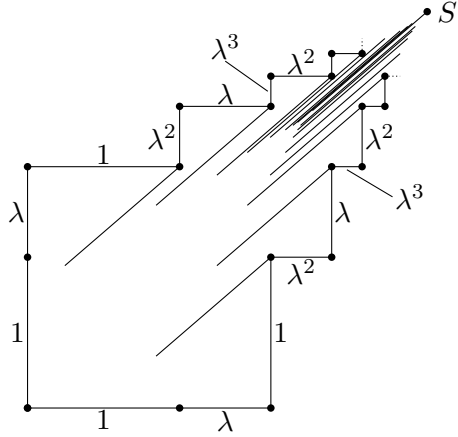
\begin{figure}
 \centering
 \begin{tikzpicture}[scale=2]
  \draw (0,0) -- node[below=-0.1]{$1$} (1,0)
   -- node[below=-0.1]{$\lambda$} (1.6,0)
   -- node[right=-0.1]{$1$} (1.6,1)
   -- node[below=-0.1]{$\lambda^2$} (2,1)
   -- node[right=-0.1]{$\lambda$} (2,1.6)
   -- (2.2, 1.6) 
   -- node[right=-0.1]{$\lambda^2$} (2.2,2)
   -- (2.35,2)
   -- (2.35,2.2);
   \draw[densely dotted] (2.35,2.2) -- (2.45,2.2);
   
  \fill (0,0) circle (0.7pt);
  \fill (1,0) circle (0.7pt);
  \fill (1.6,0) circle (0.7pt);
  \fill (1.6,1) circle (0.7pt);
  \fill (2,1) circle (0.7pt);
  \fill (2,1.6) circle (0.7pt);
  \fill (2.2,1.6) circle (0.7pt);
  \fill (2.2,2) circle (0.7pt);
  \fill (2.35,2) circle (0.7pt);
  \fill (2.35,2.2) circle (0.7pt);
   
  \draw (2.1, 1.58) -- (2.4, 1.4) node[right=-0.1]{$\lambda^3$};

  \begin{scope}[rotate=90, yscale=-1]
   \draw (0,0) -- node[left=-0.1]{$1$} (1,0)
    -- node[left=-0.1]{$\lambda$} (1.6,0)
    -- node[above=-0.1]{$1$} (1.6,1)
    -- node[left=-0.1]{$\lambda^2$} (2,1)
    -- node[above=-0.1]{$\lambda$} (2,1.6)
    -- (2.2, 1.6) 
    -- node[above=-0.1]{$\lambda^2$} (2.2,2)
    -- (2.35,2)
    -- (2.35,2.2);
   \draw[densely dotted] (2.35,2.2) -- (2.45,2.2);
    
   \fill (0,0) circle (0.7pt);
   \fill (1,0) circle (0.7pt);
   \fill (1.6,0) circle (0.7pt);
   \fill (1.6,1) circle (0.7pt);
   \fill (2,1) circle (0.7pt);
   \fill (2,1.6) circle (0.7pt);
   \fill (2.2,1.6) circle (0.7pt);
   \fill (2.2,2) circle (0.7pt);
   \fill (2.35,2) circle (0.7pt);
   \fill (2.35,2.2) circle (0.7pt);
    
   \draw (2.1, 1.58) -- (2.3, 1.3) node[above=-0.1]{$\lambda^3$};
  \end{scope}
  
  \fill (2.628, 2.628) node[right]{$S$} circle (0.7pt);
  \draw (2.628, 2.628) -- +(-139:1);
  \draw (1.6, 1) -- +(-139:1);
  \draw (2, 1.6) -- +(-139:1);
  \draw (2.2, 2) -- +(-139:1);
  \draw (2.35, 2.2) -- +(-139:1);
  \draw (2.45, 2.35) -- +(-139:1);
  \draw (2.5, 2.45) -- +(-139:1);
  \draw (2.53, 2.5) -- +(-139:1);
  \draw (2.55, 2.53) -- +(-139:1);
  \draw (1, 1.6) -- +(-139:1);
  \draw (1.6, 2) -- +(-139:1);
  \draw (2, 2.2) -- +(-139:1);
  \draw (2.2, 2.35) -- +(-139:1);
  \draw (2.35, 2.45) -- +(-139:1);
  \draw (2.45, 2.5) -- +(-139:1);
  \draw (2.5, 2.53) -- +(-139:1);
  \draw (2.53, 2.55) -- +(-139:1);
 \end{tikzpicture}
 \caption{Geodesic segments in a direction which cannot be the eigen direction of a parabolic element in the Veech group.}
 \label{fig:non_parabolic_directions}
\end{figure}

 Consider a direction $\theta$ which corresponds to a slope in $(\lambda, \lambda^{-1})$ and suppose $P$ is a parabolic element with eigen direction $\theta$. There exist infinitely many geodesic segments in direction~$\theta$ that start in the singularity and have length $1$ (see \autoref{fig:non_parabolic_directions}).
 Then $P$ has to map the set of these geodesic segments to itself. Note that the geodesic segment that starts in $S$ is distinguished from all others as it is the limit of the geodesic segments in this set. Hence, $P$ has to fix this distinguished segment pointwise.
 As $P$ is locally acting as a shear on the ladder surface, this implies that also the geodesic segments close to the limit segment have to be fixed pointwise.
 So $P$ fixes an open subset of the ladder surface. Therefore, $P$ is the identity and not a parabolic element.
\end{proof}
\end{lem}

\section{Candidate for the Veech group} \label{sec:candidate}

We found already a subgroup of the Veech group of the ladder surface with parameter $\lambda_{k,l}$: In the description by semiregular hexagons, we have the elliptic element
$\begin{pmatrix}
  \cos \sfrac{2\pi}{3} & -\sin \sfrac{2\pi}{3} \\
  \sin \sfrac{2\pi}{3} & \cos \sfrac{2\pi}{3}
 \end{pmatrix}
 =
 \begin{pmatrix}
  - \sfrac{1}{2} & - \sfrac{\sqrt{3}}{2} \\
  \sfrac{\sqrt{3}}{2} & - \sfrac{1}{2}
 \end{pmatrix}$
in the Veech group. Hence in the description as ladder surface, the corresponding Veech group element is given by
\begin{equation*}
 R \coloneqq 
 \begin{pmatrix}
  1 & -\sfrac{1}{\sqrt{3}} \\
  0 & \sfrac{2}{\sqrt{3}}
 \end{pmatrix}
 \cdot
 \begin{pmatrix}
  - \sfrac{1}{2} & - \sfrac{\sqrt{3}}{2} \\
  \sfrac{\sqrt{3}}{2} & - \sfrac{1}{2}
 \end{pmatrix}
 \cdot
 \begin{pmatrix}
  1 & -\sfrac{1}{2} \\
  0 & \sfrac{\sqrt{3}}{2}
 \end{pmatrix}
 =
 \begin{pmatrix}
  -1 & -1 \\
  1  & 0
 \end{pmatrix}
 .
\end{equation*}

Furthermore, let $T \coloneqq
 \begin{pmatrix}
  1 & k(1+\lambda) \\
  0 & 1  
\end{pmatrix}$
be the parabolic element defined by the horizontal cylinder decomposition.

Let $G = \left\langle R, T \right\rangle$ be the group generated by the elliptic and the parabolic element in~$\PSL(2,\mathbb{R})$. We determine the fundamental domain of $G$ in $\mathbb{H}$ in this section.

\bigskip

\begin{figure}[h]
 \centering
  \begin{tikzpicture}[scale=1.85]
   \fill[gray!20!white] (-2,0) -- (3.5,0) -- (3.5,1.8) -- (-2,1.8) -- (-2,0);
   \draw (-2,0) -- (-2,1.8);
   \draw (3.5,0) -- (3.5,1.8);

   \draw[->, gray] (-2.3,0) -- (3.8,0) node[right] {};
   \draw[->, gray] (0,0) -- (0,1.9) node[above] {};		

   \draw (-2,0.05) -- (-2,-0.05) node[below] {$-2$};
   \draw (-1,0.05) -- (-1,-0.05) node[below] {$-1$};		
   \draw (0,0.05) -- (0,-0.05) node[below] {$0$};
   \draw (1,0.05) -- (1,-0.05) node[below] {$1$};
   \draw (3.5,0.05) -- (3.5,-0.05) node[below] {$k(1+\lambda)-2$};
   \draw (-0.05,1) -- (0.05,1) node[right] {$\mathrm{i}$};				
   
   \draw (2,1.4) node {$F_1$};
  \end{tikzpicture}
 \caption{Fundamental domain of $\left\langle T \right\rangle$.}
 \label{fig:fundamental_domain_parabolic}
\end{figure}
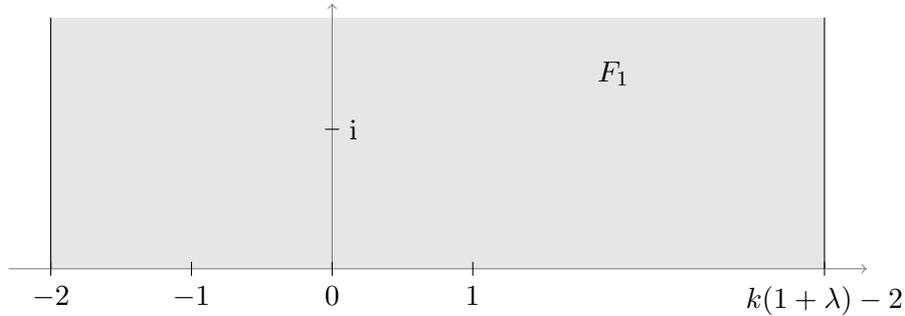

The elements of $\PSL(2,\mathbb{R})$ act as Möbius transformations on $\overline{\mathbb{H}} \coloneqq \mathbb{H} \cup\mathbb{R} \cup \infty$. Under this action, the fixed point of $T$ is $\infty$.
Hence, a fundamental domain of the subgroup $\left\langle T \right\rangle$ can be chosen to be the strip $F_1 \coloneqq \{z \in \mathbb{H} : -2 < \Re(z) < k\cdot(1+\lambda) -2\}$ as shown in~\autoref{fig:fundamental_domain_parabolic}.

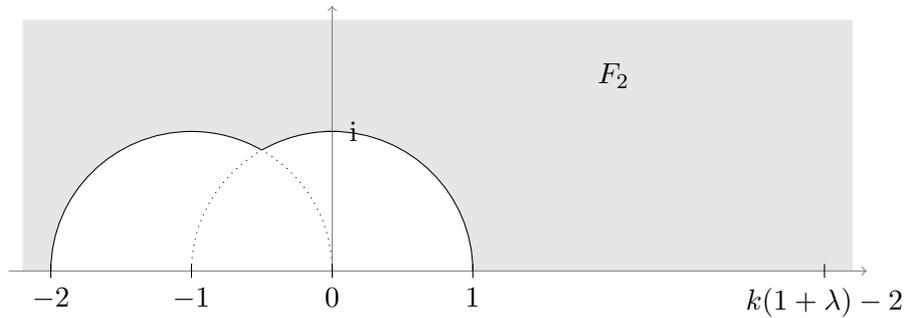
\begin{figure}[h]
 \centering
  \begin{tikzpicture}[scale=1.85]
   \fill[gray!20!white] (-2.2,0) -- (-2,0) arc (180:60:1) arc (120:0:1) -- (3.7,0) -- (3.7,1.8) -- (-2.2,1.8) -- (-2.2,0);
   \draw (-2,0) arc (180:60:1);
   \draw[dotted] (0,0) arc (0:60:1);
   \draw (1,0) arc (0:120:1);
   \draw[dotted] (-1,0) arc (180:120:1);

   \draw[->, gray] (-2.3,0) -- (3.8,0) node[right] {};
   \draw[->, gray] (0,0) -- (0,1.9) node[above] {};		

   \draw (-2,0.05) -- (-2,-0.05) node[below] {$-2$};
   \draw (-1,0.05) -- (-1,-0.05) node[below] {$-1$};		
   \draw (0,0.05) -- (0,-0.05) node[below] {$0$};
   \draw (1,0.05) -- (1,-0.05) node[below] {$1$};
   \draw (3.5,0.05) -- (3.5,-0.05) node[below] {$k(1+\lambda)-2$};
   \draw (-0.05,1) -- (0.05,1) node[right] {$\mathrm{i}$};				
   
   \draw (2,1.4) node {$F_2$};
  \end{tikzpicture}
 \caption{Fundamental domain of $\left\langle R \right\rangle$.}
 \label{fig:fundamental_domain_elliptic}
\end{figure}

For the elliptic element $R$, the fixed point in $\overline{\mathbb{H}}$ is $\frac{-1 + \sqrt{3}\cdot \mathrm{i}}{2}$. This means that $R$ sends the vertical line from $-\frac{1}{2}$ to $\infty$ to the geodesic from $1$ to $-1$. Its inverse sends the vertical line to the geodesic from $-2$ to $0$.
Hence, we have that $F_2$ as in \autoref{fig:fundamental_domain_elliptic} is a fundamental domain of the subgroup $\left\langle R \right\rangle$.

\begin{figure}
 \centering
  \begin{tikzpicture}[scale=1.85]
   \fill[gray!20!white] (-2,0) arc (180:60:1) arc (120:0:1) (3.5,0) -- (3.5,1.8) -- (-2,1.8) -- (-2,0);
   
   \draw (-2,0) -- (-2,1.8);
   \draw (3.5,0) -- (3.5,1.8);
   
   \draw (-2,0) arc (180:60:1);
   \draw[dotted] (0,0) arc (0:60:1);
   \draw (1,0) arc (0:120:1);
   \draw[dotted] (-1,0) arc (180:120:1);

   \draw[->, gray] (-2.3,0) -- (3.8,0) node[right] {};
   \draw[->, gray] (0,0) -- (0,1.9) node[above] {};		

   \draw (-2,0.05) -- (-2,-0.05) node[below] {$-2$};
   \draw (-1,0.05) -- (-1,-0.05) node[below] {$-1$};		
   \draw (0,0.05) -- (0,-0.05) node[below] {$0$};
   \draw (1,0.05) -- (1,-0.05) node[below] {$1$};
   \draw (3.5,0.05) -- (3.5,-0.05) node[below] {$k(1+\lambda)-2$};
   \draw (-0.05,1) -- (0.05,1) node[right] {$\mathrm{i}$};				
   
   \draw (2,1.4) node {$F$};
  \end{tikzpicture}
 \caption{Fundamental domain of $G$.}
 \label{fig_fundamental_domain_G}
\end{figure}
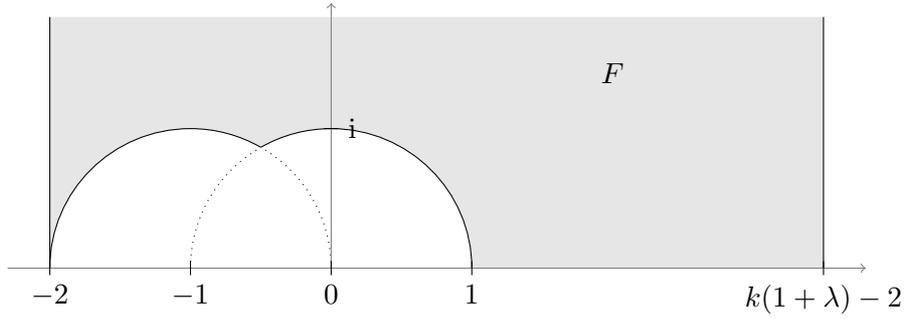

We now consider $F \coloneqq F_1 \cap F_2$ as sketched in \autoref{fig_fundamental_domain_G}.
Note that $F$ is an open convex hyperbolic polygon, $T$ and $R$ constitute a side pairing for $F$, and the cusp at the boundary of $F$, namely $\infty$, is the fixed point of a parabolic element, namely $T$.
Hence Poincaré's fundamental polyhedron theorem implies that the group generated by the side pairings, i.e.\ $G$, is discrete and that $F_G \coloneqq F$ is a fundamental domain for $G$
(see \cite{maskit_71}).

\section{The Veech group} \label{sec:Veech_group}

For translation surfaces of finite type, the Veech group is always discrete. This is not true in general for translation surfaces of infinite type but can be shown for ladder surfaces with parameter $\lambda \in (0,1)$ with the following lemma.

\begin{lem}
 Let $\Gamma$ be the Veech group of the ladder surface with parameter $\lambda \in (0,1)$.
 Then the limit set of $\Gamma$ has empty intersection with the open real segment from $\lambda$ to $\lambda^{-1}$.
 
\begin{proof}
 We first determine the slopes of eigen directions of parabolic elements that correspond to fixed points in $(\lambda, \lambda^{-1})$.
 The parabolic element
 $
  \begin{pmatrix}
   1-\lambda & 1 \\
   -\lambda^2 & 1+\lambda
  \end{pmatrix}
 \in \SL(2,\mathbb{R})
 $
 has an eigen direction with slope $\lambda$ and fixed point $\lambda^{-1} \in \partial \mathbb{H}$, and the parabolic element
 $
 \begin{pmatrix}
  1 + \lambda & -\lambda^2 \\
  1 & 1 - \lambda
 \end{pmatrix}
 \in \SL(2,\mathbb{R})
 $
 has an eigen direction with slope $\lambda^{-1}$ and fixed point $\lambda \in \partial \mathbb{H}$.
 Therefore, every parabolic element that has a fixed point in $(\lambda, \lambda^{-1})$ has an eigen direction with a slope between $\lambda$ and $\lambda^{-1}$.
 We have already seen in \autoref{lem_forbidden_eigen_directions} that no parabolic element with such an eigen direction exists.

 Suppose now that a point in $(\lambda,\lambda^{-1})$ is contained in the limit set of $\Gamma$.
 Since the fixed points of hyperbolic elements of $\Gamma$ are dense in the limit set, there also exists a point~$x \in (\lambda,\lambda^{-1})$ which is the fixed point of a hyperbolic element.
 In particular, $x$ is the attracting fixed point of a hyperbolic element $h \in \Gamma$.
 Note that $\infty$, $0$, and $-1$ are the fixed points of the parabolic elements $T$, $R^{-1} T R$, and $R T R^{-1}$, respectively, and not all three of them can be the repelling fixed point of $h$. Hence, there exists an $n\in \mathbb{N}$ such that at least one of $h^n(\infty)$, $h^n(0)$, and $h^n(-1)$ is contained in $(\lambda,\lambda^{-1})$. However, this would again correspond to a parabolic element in $\Gamma$ with a fixed point in $(\lambda,\lambda^{-1})$. Therefore, $h$ cannot be contained in $\Gamma$.
\end{proof}
\end{lem}

In particular, we can deduce the following description of $\Gamma$.

\pagebreak[3]

\begin{prop}[Veech group as Fuchsian group]
 The Veech group $\Gamma$ of the ladder surface with parameter $\lambda \in (0,1)$ is a non-elementary Fuchsian group of the second kind.

\begin{proof}
 As $\Gamma$ acts discontinuously on a neighborhood of $(\lambda, \lambda^{-1})$ in $\mathbb{H}$, the group $\Gamma$ is a discrete subgroup of $\PSL(2,\mathbb{R})$, i.e.\ a Fuchsian group.
 
 On one hand, the fixed points of the parabolic elements $T$, $R^{-1} T R$, and $RTR^{-1}$ are different elements of the limit set. This implies that the limit set contains infinitely many elements and $\Gamma$ is not elementary.
 On the other hand, the limit set of $\Gamma$ is disjoint from~$(\lambda,\lambda^{-1})$, hence $\Gamma$ is of the second kind. 
\end{proof}
\end{prop}

We need now one more lemma which implies that the candidate group $G$ from the previous section is a normal subgroup of the Veech group $\Gamma$. For this, we restrict to~$\lambda = \lambda_{k,l}$ with $l=1$. Note that this also implies $\lambda \in (0,1)$.

\begin{lem}
 Let $\lambda = \lambda_{k,1}$.
 For any $A\in \Gamma$, the conjugates $A^{-1} T A$ and $A^{-1} R A$ are contained in~$G$.
 
\begin{proof}
 Let $A\in \Gamma$. Then $A^{-1} T A \in \Gamma$ is parabolic and has a fixed point $c \in \mathbb{R} \cup \{\infty\}$. There exists an element $g \in G$ such that $g \cdot c$ is contained in the closure of $F_G$.
 Note that $F_G$ has exactly one cusp and a hole which cannot contain fixed points of parabolic elements of $\Gamma$, so $g \cdot c$ has to be $\infty$. In contrast, $g \cdot c$ is also the fixed point of the parabolic element $g A^{-1} T A g^{-1}$ which implies that $g A^{-1} T A g^{-1}$ has the same eigen direction as $T$.
 Hence, $g A^{-1} T A g^{-1}$ has to fix the horizontal cylinder decomposition from \autoref{fig:horizontal_cylinder_decomposition}. In this decomposition, there exists a widest horizontal cylinder with circumference $1 + \lambda + \lambda^2$. Its boundary consists of four saddle connections of lengths $1$, $\lambda + \lambda^2$, $1 + \lambda$, and $\lambda^2$. Hence, this boundary has to be fixed pointwise which means that every parabolic element with horizontal eigen direction has to be a power of $T$.
 In particular, $g A^{-1} T A g^{-1}$ is a power of~$T$. Hence, $g A^{-1} T A g^{-1}$ is contained in $G$ and so is $A^{-1} T A$.

 Furthermore, $A^{-1} R A$ is elliptic of order $3$ and has a fixed point $c \in \mathbb{H}$. There exists an element $g \in G$ such that $g \cdot c$ is either the fixed point of $R$ or is contained in the interior of~$F_G$.
 The image of the geodesics from the fixed point of $R$ to $\infty$, $-1$, and $0$ under $g A^{-1}$ are geodesics from $g \cdot c$ to cusps. Such geodesics cannot exist for a point in the interior of~$F_G$, hence $g \cdot c$ has to be the fixed point of $R$. Therefore, $R$ and $g A^{-1} R A g^{-1}$ are both elliptic elements of the same order and with the same fixed point. Hence, they generate the same subgroup and therefore $g A^{-1} R A g^{-1}$ is contained in $G$ and so is $A^{-1} R A$.
\end{proof}
\end{lem}

\begin{thm}[Generators for Veech group]
 The Veech group $\Gamma$ of the ladder surface with parameter $\lambda_{k,1}$ is generated by the parabolic element
 $P =
 \begin{pmatrix}
  1 & k\cdot (\lambda+1) \\
  0 & 1
 \end{pmatrix}$ and the elliptic element
 $R =
 \begin{pmatrix}
  -1 & -1 \\
   1 &  0
 \end{pmatrix}$.

\begin{proof}
 Let $F_\Gamma$ be a fundamental domain of $\Gamma$ with $F_\Gamma \subseteq F_G$ and such that $F_G$ is tessellated by translates of $F_\Gamma$ under elements of $\Gamma$.
 As $F_G$ contains exactly one cusp and a hole and as $\Gamma$ cannot have cusps in the hole of $F_G$, also $F_\Gamma$ needs to have one cusp, which is $\infty$.

 Let $A\in \Gamma$. Then we can choose $g \in G$ such that $g A^{-1} \cdot F_{\Gamma} \subseteq F_G$. As in the previous lemma, $g A^{-1} T A g^{-1}$ is a power of $T$ and $g A^{-1} R A g^{-1}$ is a power of $R$.

 Note that, if for an elliptic or parabolic element $H \in \PSL(2,\mathbb{R})$ and $n,m \in \mathbb{Z}$, we have $g A^{-1} H^n A g^{-1} = H^m$, then $gA^{-1}$ has the same fixed point as $H$. Hence, $g A^{-1}$ and $H$ are both powers of the same element and hence they commute. This shows that $n=m$.

 In particular, we have $g A^{-1} T A g^{-1} = T$ and $g A^{-1} R A g^{-1} = R$, hence $g A^{-1}$ commutes with both $T$ and $R$. Recall that nontrivial Möbius transformations commute if and only if they have the same set of fixed points. Therefore, we conclude that $g A^{-1}$ is the identity and hence $A \in G$.
\end{proof}
\end{thm}

If we consider the Veech group in $\SL(2,\mathbb{R})$ instead of in $\PSL(2,\mathbb{R})$, also the rotation
$\begin{pmatrix}
  -1 & 0 \\
  0 & -1
\end{pmatrix}$
might be contained in the Veech group.
This element would have to preserve the horizontal saddle connection of maximal length $1+\lambda$ and the vertical saddle connection of maximal length $1+\lambda$. However, the intersection pattern of these saddle connections is not preserved under a rotation by $\pi$, hence 
$\begin{pmatrix}
  -1 & 0 \\
  0 & -1
\end{pmatrix}$
cannot be an element of the Veech group.

In contrast,
$\begin{pmatrix}
  0 & 1 \\
  1 & 0
\end{pmatrix}$
is an affine orientation-reversing homeomorphism of the ladder surface.
Hence, the extended Veech group $\Gamma^{\pm}$ is generated by $P$, $R$ and
$\begin{pmatrix}
  0 & 1 \\
  1 & 0
\end{pmatrix}$.

\paragraph{Acknowledgements} I would like to thank Pat Hooper for asking and discussing the question.

\bibliographystyle{amsalpha}
\bibliography{/home/anja/Documents/Mathematik/literature/BibTex/Literatur}

\end{document}